\documentclass[12pt]{amsart}
\newtheorem{Thm}{Theorem}
\newtheorem{thm}{Theorem}[section]
\newtheorem{lem}[thm]{Lemma}

\theoremstyle{definition}
\newtheorem{defn}{Definition}
\newtheorem{defnn}{Definition}

\newcommand{\R}{\mathbb R}

\newcommand{\Z}{\mathbb Z}

\newcommand{\sinc}{\mbox{sinc}}
\newcommand{\tl}{\tau_\lambda\, }
\newcommand{\shat}{{}^{\wedge}}

\begin{document}

\title[]{Completeness in $L^1 (\R)$ of discrete translates}%
\author{Joaquim Bruna}
\thanks{The first author is supported by projects BFM2002-04072-C02-02 and
2001SGR00172}
\author{Alexander Olevskii}
\thanks{The second author is partially supported by a grant of Israel Scientific Foundation}
\author{Alexander Ulanovskii}

\address{Departament de Matem\`{a}tiques, Universitat Aut\`{o}noma de
Barcelona\\ 08193 Cerdanyola del Vall\`{e}s, SPAIN\\ bruna@mat.uab.es}
\address{School of Mathematical Sciences, Tel-Aviv University \\ Ramat Aviv 69978, Israel\\olevskii@post.tau.ac.il}
\address{Stanavanger University College\\P.O.Boks 2557 Ullandhaug\\Stavanger 4091, Norway\\alexander.ulanovskii@tn.his.no}


\begin{abstract} We characterize, in terms of the Beurling-Malliavin density, the discrete
spectra $\Lambda\subset\R$ for which a generator exists, that is a
function $\varphi\in L^1(\R)$ such that its $\Lambda$-translates
$\varphi(x-\lambda), \lambda\in\Lambda$, span $L^1(\R)$. It is
shown that these spectra coincide with the uniqueness sets for
certain analytic classes. We also present examples of discrete
spectra  $\Lambda\subset\R$ which do not admit a single generator
while they admit a pair of generators.
\end{abstract}

\maketitle
\section{Introduction}

\noindent 1.1. The famous Wiener Tauberian theorem states that a
function $f\in L^1 (\R)$ spans $L^1(\R)$ by translations, in the
sense that the linear combinations of translates $(\tl
f)(t)=f(t-\lambda)$ of $f$, $\lambda \in \R$, are dense in $L^1
(\R)$, if and only if $\hat f(\zeta)\neq 0$ for all $\zeta \in
\R$. The corresponding result for $L^2(\R)$ is that a function
$f\in L^2(\R)$ spans $L^2(\R)$ by translations if and only if
$\hat f (\zeta)\neq 0$ almost everywhere in $\R$.

 Let $X$ be some translation-invariant function space on $\R$
(that is, $\tl f \in X$ for $f\in X$, $\lambda \in \R$). A
function $\varphi\in X$ may have the property that only a certain
set of translates $\tau_\lambda \varphi$ where $\lambda$ belong to
some set $\Lambda,$ suffice to span $X$:

\begin{defn} Let $\varphi \in X$ and $\Lambda\subset\R$.
We say that \emph{$\varphi$ is a $\Lambda$-generator for $X$} if
the linear span $T(\varphi,\Lambda)$ of the translates $\tl
\varphi$, with $\lambda \in \Lambda$, is dense in $X$. \end{defn}

It is natural to ask which spectra $\Lambda$ admit a generator
$\varphi$ in a fixed space $X$. This question is most interesting
when $\Lambda$ is discrete, which we will assume from now on.

In case $\Lambda$ is the set of integers $\Z$, that is, we are
dealing with \emph{integer translates} of a fixed function, it is
well-known that no $\Z$-generators exist in $L^p(\R), 1\leq p\leq
2$. In $L^2(\R)$, this easily follows from Plancherel's theorem
and the fact that $T(\varphi, \Lambda)$ has Fourier transform
\begin{equation}\label{1}
  T(\varphi, \Lambda)\shat {}
  = \hat\varphi \;\mathcal{E} (\Lambda)
\end{equation}
where $\mathcal{E} (\Lambda)$ denotes the linear span of the
exponentials $e^{i \lambda \zeta}$ with frequencies $\lambda \in
\Lambda$. If $\Lambda =\Z$ this consists entirely of
$2\pi$-periodic functions and hence $\hat\varphi\cdot\mathcal{E}
(\Lambda)$ cannot be dense in $L^2(\R)$. A similar argument works
in $L^p(\R), 1\leq p<2$.

However, surprisingly enough, in $L^p(\R)$, $p>2$, there do exist
$\Z$-generators. This result was established in \cite{AO}, and
another proof can be obtained from results of \cite{N} (see also
\cite{F} for a particular case).

In the space $L^2(\R)$, using a certain construction based on
small divisors, Olevskii \cite{O} showed that an arbitrary
perturbation of $\Z$ of the form
\begin{equation}\label{2}
  \Lambda = \{\, n+a_n\; , \; a_n \neq 0\; , a_n \rightarrow 0\,\}
\end{equation}
admits a generator $\varphi \in L^2 (\R)$.

It is immediate to see using Plancherel's theorem that $T(\varphi,
\Lambda)$ is dense in $L^2(\R)$ if and only if $\hat
\varphi(\zeta)\neq 0$ almost everywhere and  $\mathcal{E}
(\Lambda)$ is dense in the weighted $L^2$-space $L^2(\R,\omega)$,
with $\omega=|\hat\varphi |^2$, an a.e. positive weight. In
particular, if $E_{\varepsilon}$ denotes the set
$E_{\varepsilon}=\{\omega\geq \varepsilon\}$, $\mathcal{E}
(\Lambda)$ will be dense in $L^2 (E_{\varepsilon})$, and
$|E_{\varepsilon}|\rightarrow \infty$ as $\varepsilon \to 0$.
Thus, if $\Lambda$ has a generator in $L^2(\R)$, $\mathcal{E}
(\Lambda)$ is dense in $L^2$ in sets of arbitrarily large measure.

This shows the connection of these questions with the subject of
\emph{density of exponentials} $\mathcal{E} (\Lambda)$ in function
spaces and, in particular, with Landau's results. Landau \cite{La}
constructed certain perturbations of the integers
$\Lambda=\{n+a_n\}$ where $a_n$ are bounded, such that
$\mathcal{E} (\Lambda)$ is dense in $L^2$ on any finite union of
the intervals $(2\pi(k-1)+\varepsilon, 2\pi k -\varepsilon)$,
$\varepsilon >0$, in particular on sets with arbitrarily large
measure. In \cite{U} Landau's result was extended to every
sequence $\Lambda$ as in (2), where $a_n$ have an exponential
decay. We mention here that if $\mathcal{E} (\Lambda)$ is complete
in $L^2$ on `Landau sets', then one can construct a
$\Lambda-$generator for $L^2(\R)$ which belongs to the Schwartz
class  $S(\R)$. Such generators are presented in \cite{OU} for
sequences (2) with exponentially small $a_n$. It is also shown in
\cite{OU} that the exponential decay is in a sense necessary,
since a slower decay of $a_n$ cannot guarantee existence of a
generator even from $L^1(\R).$

In general there is a sort of balance between the size of
$\Lambda$ and the `smallness' of $\hat\varphi$. The faster
$\hat\varphi$ tends to zero at $\pm\infty$ the sparser spectra
$\Lambda$ may serve as translation sets, the "denser" $\Lambda$ is
the more general $\varphi$ may work as a generator. The spectra
$\Lambda$ considered in \cite{O} and \cite{OU} are "sparse" in the
sense that they all have density one; a number of results for
spectra $\Lambda$ with infinite density can be found in \cite {Z}
and \cite{S}.

In connection with these questions, we point out that no Riesz
bases (nor a frame) exists in $L^2(\R)$ consisting of translates
of a fixed function $\varphi$ (\cite{CDH}). On the other hand, to
the best of our knowledge, it is not known whether a Schauder
bases of translates exists in $L^2(\R)$.

\vspace{0.5cm}

\noindent 1.2. This paper deals with the case $X=L^1(\R)$. Our
main result (Theorem 1 below) gives a characterization of the
translation sets $\Lambda \subset \R$ admitting a generator
$\varphi$. The $L^1$-case is in a sense easier than the $L^2$-case
because now $\hat\varphi$ must be a non-vanishing continuous
function and what is involved is the question of density of
exponentials $\mathcal{E} (\Lambda)$ in \emph{intervals}. It is
therefore not surprising that the \emph{Beurling-Malliavin
spectral radious formula}, which we now recall, appears in this
setting and in the statement of the main result.

\noindent The spectral radius $R(\Lambda)$ of a set $\Lambda
\subset \R$ is defined
$$
R(\Lambda)=\sup \{ \,\rho >0 \,; \; \mathcal{E} (\Lambda) \text{
is complete in} \, C[-\rho, \rho\,]\,\},
$$
where $C(I)$ denotes the space of continuous functions on the
interval $I$. One sets $R(\Lambda)=0$ when $\mathcal{E}(\Lambda)$
is not complete in $C[-\rho, \rho\,]$ for any positive $\rho$.
Thus if $0<R(\Lambda)<\infty,$ we have that $\mathcal{E}
(\Lambda)$ is complete in $C[-\rho, \rho]$ if $\rho < R(\Lambda)$
and incomplete if $\rho
> R(\Lambda)$. The Fourier transform of the dual space, the space
of finite complex Borel measures supported in $[-\rho, \rho]$, is
a (proper) subspace of the Bernstein space
$$
B_\rho = \bigl\{\, F \text{ entire: $ |F(x+iy)|\leq C_F e^{\rho
|y|}$}, \, x+iy\in\mathbb{C} \bigr\},\,
$$
where $C_F>0$ is a constant depending on $F$. If $\Lambda$ is not
complete in $C[-\rho,\rho]$, then there is a function $F\in
B_\rho$ which vanishes on $\Lambda$. On the other hand, if a
function $G\in B_\rho$ vanishes on $\Lambda$, the function
$F(t)=G(t)(\sin \varepsilon t/\varepsilon t)^2$ is the Fourier
transform of a finite (absolute continuous) measure concentrated
on $[-\rho-2\varepsilon,\rho+2\varepsilon]$ which is orthogonal to
$\mathcal{E} (\Lambda)$. Hence, for $R(\Lambda)$ finite and
positive, we have
\newcommand{\conj}[2]{\bigl\{\, {#1} \,:\,{#2}\,\bigr\}}
\begin{align}\label{3}
R(\Lambda)&=\sup\,
    \conj{\rho >0}{\text{$\Lambda$ is a uniqueness set for $B_\rho$}}
                  = \\[3pt] \nonumber
 &=\inf\, \conj{\rho >0}{
\text{$\exists F\in B_\rho, F\ne 0,$ such that $F_{|\Lambda} =0
$}}\,.
\end{align}
Here and in the sequel, $\Lambda$ being a \emph{uniqueness set for
a certain class $Y$} means that $F\in Y$ and $F_{|\Lambda}=0$
imply $F\equiv 0$.

In the beginning of the sixties  Beurling and Malliavin (see
\cite[vol. 2]{Ko}) established that $R(\Lambda)=\pi
D_{BM}(\Lambda)$, where $D_{BM} (\Lambda)$ is a certain exterior
density whose definition will be recalled in section~3.

\vspace{0.5cm}

\noindent 1.3. With this notation we can now formulate our main
result:
\begin{Thm}\label{A}
A discrete set $\Lambda \subset \R$ admits a generator for
$L^1(\R)$ if and only if $R(\Lambda)=+\infty$.
\end{Thm}

The necessity of the condition $R(\Lambda)=+\infty$ follows
essentially from the remarks in paragraph 1.1. and will be
detailed in section 2. We will give two proofs of the sufficiency.
The first one, also in section 2, is constructive, based solely on
$R(\Lambda)=+\infty$ without using any density. The second one, in
section 3, uses the geometrical information
$D_{BM}(\Lambda)=+\infty$.

This second proof develops the following natural approach to the
problem, relating it with uniqueness sets for \emph{generalized
Bernstein classes}. Let $N$ be the class of all functions $f\in
L^1(\R)$ such that $\hat f$ does not vanish:
\begin{equation}\label{N}
N=\{f\in L^1(\R): \hat f(\zeta)\ne 0, \zeta\in\R\}.
\end{equation}
By duality, a function $\varphi$ is a $\Lambda$-generator for
$L^1(\R)$ if and only if $h\in L^\infty (\R)$ and $(h \star
\check{\varphi})(\lambda)= \int_{-\infty}^{+\infty} h(t)\varphi
(t-\lambda)\, dt=0$, $\lambda \in \Lambda,$ implies $h=0$
  (here $\check{\varphi} (t)=\varphi (-t) $). Applying Wiener's theorem,
we see that  $\varphi$ is a $\Lambda$-generator   if and only if
$\varphi\in N$, and  $h\in L^\infty (\R)$ and $(h \star
\check{\varphi})(\lambda)=0, \lambda \in \Lambda,$ imply
$h\star\check\varphi=0$.  This can be restated by saying that
\emph{$\varphi$ is a $\Lambda$-generator in $L^1(\R)$ if and only
if it satisfies: (i) $\varphi\in N$, and (ii) $\Lambda$ is a
uniqueness set for the class $L^\infty (\R) \star \check
\varphi$.}

Given any $\Lambda$ with infinite spectral radius, it follows from
(\ref{3}) that $\Lambda$ is a uniqueness set for every Bernstein
class $B_\rho$. However, clearly, we have $N\cap
B_\rho=\emptyset,$ for all $\rho$, and so no space $L^\infty (\R)
\star \check \varphi, \varphi\in N,$ is included in any of the
Bernstein spaces. In order to prove Theorem 1, one may try to
construct \emph{larger} spaces of analytic functions $Y$ such that
$\Lambda$ is a uniqueness set for $Y$, and $L^\infty (\R) \star
\check \varphi\subseteq Y$, for some $\varphi\in N$. It turns out
that the `smallest' spaces $Y$ with this property can be realized
as the generalized Bernstein classes $B_{\sigma}$ defined as
follows. Let $\sigma$ be a monotone function on $(0,\infty)$
 satisfying $\sigma (y)\nearrow\infty, y\to\infty$. We set:
\begin{equation}\label{bs}
B_{\sigma}=\conj{F\text{ entire}}{|F(x+iy)|\leq C_F\, e^{|y|\sigma
(|y|)}, x+iy\in \mathbb{C}\,}.
\end{equation}

In section 3 we prove:

\begin{Thm}\label{B}
For a discrete set $\Lambda \subset \R$, the following conditions
are equivalent:
\begin{itemize}
  \item[(a)] $\Lambda$ admits a generator $\varphi$ for $L^1(\R)$.
  \item[(b)] $\Lambda$ is a uniqueness set for some generalized
  Bernstein class $B_\sigma$.
\end{itemize}
When \emph{(b)} holds, the generator $\varphi$ can be chosen in
$B_\sigma$.
\end{Thm}

If a set $\Lambda$ has finite spectral radius then, by Theorem 1,
$\Lambda$ does not admit generators in $L^1(\R)$. In section 4 we
look at the problem of \emph{two} generators and show that there
exist spectra $\Lambda$ which do not admit a single generator
while they admit a pair of generators. More specifically, we show
that every $\Lambda$ of form (\ref{2}) with `exponentially small'
$a_n$ admits two generators.

\section{The spectral radius proof of the main theorem}

A basic well know remark concerning the spectral or completeness
radius $R(\Lambda)$ is that its value is not affected if the
sup-norm on $[-\rho, \rho]$ is replaced by any reasonable norm.
This fact is behind both the proof of the necessity and
sufficiency of the condition.

\vspace{0.5cm}

\noindent  2.1. Let us prove first that $R(\Lambda)=+\infty$ is
necessary. If $\varphi \in L^1(\R)$ is a $\Lambda$-generator, it
follows from (1) and the trivial estimate $\Vert \hat f
\Vert_\infty \leq \Vert f\Vert_1$ that an arbitrary $\hat f$, with
$f\in L^1(\R)$, can be approximated in the sup-norm by functions
in $\hat \varphi\, \mathcal{E} (\Lambda)$. Since $\varphi \star f
\in L^1(\R)$, $\hat\varphi \hat f$ can be approximated as well. We
know too that $\hat \varphi$ has no zeros, by Wiener's theorem.
Now fix $\rho >0$; every test function $\psi$ supported in
$(-\rho, \rho)$ serves as $\hat f$, and therefore $\hat \varphi
\psi$ is well approximated by $\hat \varphi \,\mathcal{E}
(\Lambda)$ in the sup-norm. Since $\hat\varphi$ is bounded below
on $(-\rho, \rho)$ it follows that every $\psi$ is approximated in
the sup-norm by $\mathcal{E} (\Lambda)$. The density of such
$\psi$ in $C\, [-\rho, \rho]$ shows that $\mathcal{E} (\Lambda)$
is dense in $C\, [-\rho, \rho]$ and $\rho$ being arbitrary, one
has $R\,(\Lambda)=+\infty$.

\vspace{0.5cm}

\noindent 2.2. To prove that $R\,(\Lambda)=+\infty$ is a
sufficient condition we will make use of the remark above and
consider instead of the sup-norm the Sobolev norm
$$
\parallel \psi \parallel_2 + \parallel\psi'\parallel_2.
$$
More precisely, if $I$ is an interval we consider the space
$$
L_1^2 (I)= \, \conj{h\in L^2 (I)}{h' \in L^2(I)}
$$
(consisting of absolutely continuous functions) endowed with the
norm \linebreak  $\parallel h \parallel_I \, = \,
\parallel h \parallel_{L^2(I)} + \parallel h'\parallel_{L^2(I)}$.
The reason to consider $\parallel h\parallel_I$ is that if
$\parallel\hat f\parallel_{\R}$ is finite then $f\in L^1 (\R)$ and
\begin{equation}\label{6}
\parallel f \parallel_1 \, \leq \,  \parallel \hat f
\parallel_{\R},
\end{equation}
as it is immediately checked.

The fact that if $R\,(\Lambda)=+\infty$, $\mathcal{E} \,(\Lambda)$
is dense in $L_1^2\,(-\rho, \rho)$ for every $\rho >0$ can be
checked in an elementary way, i.e. by first approximating $h'$ and
then integrating.

We consider now the space
$$
E=\conj{f\in L^1(\R)}{\hat f \text{ is of class $C^1$ and
compactly supported}}.
$$
It is clear that $E$ is dense in $L^1(\R)$ and separable. We fix a
sequence $(f_k)_{k=1}^\infty$, $f_k \in E$ dense in $L^1(\R)$. Fix
some sequence $\varepsilon_k \to 0$, say $\varepsilon_k =2^{-k}$.
Associated to $\Lambda$ with $R\,(\Lambda)=+\infty$ and to the
sequences $(f_k)_{k=1}^\infty$, $(\varepsilon_k)_{k=1}^\infty$ we
will construct a positive $\Phi\in L_1^2\,(\R)$ with the following
property:

``For each $k=2, 3, ...$ there exists a trigonometric polynomial
$P_k (\zeta)=\sum_{\lambda \in \Lambda_k} c_\lambda e^{i\lambda
\zeta}$ with frequencies in a finite subset $\Lambda_k$ of
$\Lambda$, such that $\parallel \hat f_k - \, P_k
\Phi\parallel_{\R}< \varepsilon_k$".

\vspace{0.5cm}

\noindent If $\hat \varphi =\Phi$ then (6) gives $\varphi \in L^1
(\R)$ such that
$$
\parallel f_k - \sum_{\lambda \in \Lambda_k} c_\lambda \tau_\lambda
 \varphi \parallel_1 \, \leq \,
\varepsilon_k.
$$
Since given an arbitrary $g\in L^1(\R)$ and $\varepsilon >0$ there
are infinitely many $k$ such that $\parallel g - f_k
\parallel_1 < \varepsilon$, this will prove that $\varphi$ is
indeed a $\Lambda$-generator in $L^1(\R)$.

We will use at certain points the trivial estimate
$$
\parallel HG \parallel_I \; \leq \, B(G) \parallel H \parallel_I,
\ \ G, H\in L^2_1(\R),
$$
where $B(G)=\parallel G \parallel_\infty + \parallel G'
\parallel_\infty$.

Let $(I_k)_{k=1}^\infty$, $(J_k)_{k=1}^\infty$ be open intervals
centered at 0, with
$$
J_1 \subset I_1 \subset J_2 \subset I_2 \subset \ldots \subset J_k
\subset I_k \subset J_{k+1} \subset I_{k+1} \subset \ldots
$$
$I_k\setminus J_k$ consisting in two unit intervals and such that
$\hat f_k$ is supported in $J_k$. The function $\Phi$ will be an
even continuous piecewise linear function that will be built step
by step (the condition of piecewise continuity is not essential
for the construction). We will exploit the fact that
\begin{equation}\label{7}
\text{`` $\mathcal{E}  \,(\Lambda)$ is dense in $L_1^2 (I)$ for
every interval $I$ "}
\end{equation}

Let $\delta_k >0$ be such that $\sum_{n=k}^\infty \delta_n =
\varepsilon_k$. We shall now construct even piecewise linear
functions $G_k,\,k= 1, 2, ...,$ and trigonometrical polynomials
$P_k$ with frequencies from $\Lambda$, each $G_k$ being positive
and decreasing on $I_{k}\cap (0,\infty),\, G_k=0$ outside $I_{k}$,
$B(G_k)<1$, $\parallel G_k\parallel_\R\leq 1$, $G_k=G_{k+1}$ in
$I_{k-1}, k= 2, 3, ...$ and such that
\begin{equation}\label{eq1}
\parallel \hat f_k-P_k G_k\parallel_\R\leq \delta_k, \ k= 1 ,2,
...,
\end{equation}
and
\begin{equation}\label{eq2}
\max\{1, B(P_1),...,B(P_k)\}\parallel  G_{k+1}-G_k
\parallel_\R\leq \delta_{k+1}, \ k= 1, 2, ...
\end{equation}

To begin, clearly there exists an even piecewise linear function
$G_1$ positive, decreasing on $I_1\cap(0,\infty)$ and vanishing
outside $I_1$ such that $B(G_1)< 1 $ and $\parallel
G_1\parallel_\R<1.$
 Recall that $\hat f_1$ is zero outside $J_1\subset I_1$.
 By (\ref{7}),  there is a trigonometric polynomial $P_1$
with frequencies from $\Lambda$ such that $\parallel \hat
f_1/G_{1}-P_1\parallel_{I_1}<\delta_1$. Since $G_{1}$ vanishes
outside $I_1$ and $B(G_1)< 1$, this implies (\ref{eq1}) with
$k=1$.

 Next, one can choose an even piecewise linear function $G_2$ positive decreasing on
$I_2\cap(0,\infty)$ and vanishing outside $I_2$ such that $G_2$ is
so close to $G_1$ (hence, $G_2$ is close to zero on $I_2\setminus
I_1$) that $B(G_2)< 1 $, $\parallel G_2\parallel_\R<1$ and
(\ref{eq2}) holds with $k=1$. The same argument we used for $G_1$
shows that there exists $P_2$ such that (\ref{eq1}) holds with
$k=1$.

Assume that we have found $G_1,...,G_{n}$ satisfying the
properties above. Then, clearly, we may find $G_{n+1}$ which is
even, piecewise linear, positive and decreasing on
$I_{n+1}\cap(0,\infty)$ and vanishing outside $I_{n+1}$ such that
$G_{n+1}=G_{n}$ on $I_{n-1}$ and $G_{n+1}$ is so close to $G_{n}$
that $B(G_{n+1})< 1 $, $\parallel G_{n+1}\parallel_\R<1$ and
(\ref{eq2}) holds with $k=n$. Again, the argument we used for
$G_1$ shows that there is a trigonometrical polynomial $P_{n+1}$
with frequencies from $\Lambda$ such that (\ref{eq1}) holds true
with $k=n+1$.

Now define $\Phi=G_{k}$ on $I_{k-1}$. Then, $$\parallel
\Phi\parallel_\R \leq \lim_{k\to\infty}\parallel
G_k\parallel_\R\leq 1,$$ which shows that $\varphi\in L^1.$
Moreover, by (\ref{eq1}) and (\ref{eq2}), since $\hat f_k=0$
outside $J_k\subset I_k$ and $G_n=0$ outside $I_n$, we have for
every $k= 2, 3, ...$
$$\parallel \hat f_k-P_k\Phi\parallel_{\R}=
\parallel \hat f_k-P_kG_{k+1}\parallel_{I_{k}}+\sum_{n=k}^\infty \parallel
P_kG_{n+2}\parallel_{I_{n+1}\setminus I_{n}}\leq$$
$$\parallel \hat f_k-P_kG_k\parallel_{\R}+\parallel  P_kG_{k+1}-P_k G_k\parallel_{\R}+
\sum_{n=k}^\infty \parallel P_kG_{n+2}-P_kG_{n+1}\parallel_{\R}$$
$$ \leq \delta_{k}+
 \delta_{k+1}+\sum_{n=k+1}^\infty \delta_{n+1}=\varepsilon_k,
$$
as desired.

\newpage

\section{Uniqueness sets for generalized Bernstein classes}

In this section we prove Theorem 2 stated in the Introduction.

\vspace{0.5cm}

\noindent 3.1. We first verify the `easy' part  (b) $\Rightarrow$
(a). We assume that $\Lambda$ is a uniqueness set for some
$B_\sigma$ with $\sigma(y)\nearrow\infty$ as $y \to + \infty$, and
we want to produce $\varphi \in B_\sigma$ which is a
$\Lambda$-generator. The proof is based on the fact that for any
$\sigma(y)\nearrow\infty$ we have $B_\sigma\cap N\ne\emptyset$,
where the class $N$ is defined in (\ref{N}). More precisely, we
need the following

\begin{lem}\label{l00}
For any $\sigma(y)\nearrow\infty, y\to\infty,$ there is a positive
function on $(0,\infty)$,$\omega(s)\nearrow\infty$,  such that
$$\int_0^\infty e^{ys-\omega(s)}\, ds\leq e^{y\sigma(y)}, \ y\geq
0.$$
\end{lem}

This  is a simple and well-known fact, and so we omit the proof.

The proof of (a) $\Rightarrow$ (b) is as follows. First we apply
Lemma \ref{l00} to $\sigma(y)-2\epsilon$, with  some fixed number
$\epsilon>0$. Define $g$ as $\hat g(s)=e^{-\omega(|s|)}$, then
$g\in B_{\sigma - 2\epsilon}$. Set
$\varphi(t)=g(t)\,\sinc^2(\epsilon t),$ where $\sinc\, t=\sin
t/t.$ The Fourier transform of $\sinc\, \epsilon t$ is
$\chi_\epsilon(s)$, the characteristic function of the interval
$[-\epsilon,\epsilon]$. Hence, $\hat\varphi=\hat
g\star\chi_\epsilon\star\chi_\epsilon.$ Since $\hat g$ is
everywhere positive, the same is true for $\hat\varphi$, which
gives $\varphi\in N$.

We shall use the trivial estimate:
\begin{equation}\label{sinc}\left|\sinc^2\epsilon(x+iy)\right|
\leq C\frac{e^{2\epsilon|y|}}{1+x^2+y^2},\
x+iy\in\mathbb{C},\end{equation} where $C$ is some constant. This
estimate shows that $\varphi\in B_\sigma$ and that $\varphi(x+iy)$
is in $L^1$ on any line $z=x+iy, -\infty<x<\infty$. The latter
implies that for every $f\in L^{\infty}(\R)$ the convolution
$f\star\check\varphi$ is defined at every complex point $x+iy$.
Using (\ref{sinc}) and $g\in B_{\sigma-2\epsilon}$ we obtain:
$$|f\star\check\varphi(x+iy)|=\left|\int_{-\infty}^\infty \varphi(x+iy-s)f(s)\,ds\right|\leq \Vert f\Vert_\infty\int_{-\infty}^\infty|\varphi(x+iy)|\, dx\leq$$
$$\Vert f\Vert_\infty\max_x|g(x+iy)|\int_{-\infty}^\infty  C\frac{e^{2\epsilon|y|}}{1+x^2+y^2}\,
dx\leq C_1 e^{|y|\sigma(|y|)}, \exists C_1>0, \
x+iy\in\mathbb{C}.$$ Hence, $L^\infty(\R)\star \varphi\subseteq
B_\sigma$. As explained in paragraph 1.3 this shows that $\varphi$
is a $\Lambda-$generator for $L^1(\R)$.

\vspace{0.5cm}

\noindent 3.2. For the proof of $\text{(a) $\Rightarrow$ (b)}$ we
start recalling the definition of the Beurling-Malliavin exterior
density. There are several equivalent definitions of this density
$D_{BM}$. One suitable for us here is as follows (see~\cite[vol.
2]{Ko}). Suppose $\Lambda \subset (0, \infty)$ and let $D$ be a
positive number. A family of disjoint intervals $I_k=(a_k, b_k)$,
$0 < a_1 < b_1 < \cdots , a_k \nearrow + \infty$ is called
\emph{substantial for $D$} if
$$
\newcommand{\igdef}{\mathop{\mathgroup\symoperators =}}
\frac{n_{\Lambda}(a_k,b_k)}{b_k-a_k}> D\; ,\ k=1, 2, \ldots, \quad \sum_k \bigl (
\frac{b_k -a_k}{b_k}\bigr )^2 = + \infty \,.
$$
Here $n_\Lambda (I)$ is the number of points of $\Lambda$ in an
interval $I$. The density $D_{BM}(\Lambda)$ is then defined,
$$
D_{BM} (\Lambda)= \sup \, \conj{D>0}{\text{there exists  a
substantial family for $D$}}\,.
$$
If no $D>0$ admits a substantial family, one sets $D_{BM}
(\Lambda)=0$. For a general $\Lambda$, one defines
$D_{BM}(\Lambda)=\max \{ D_{BM}(\Lambda^+),\, D_{BM}(\Lambda^-)\}$
where $\Lambda^+ = \Lambda \cap \R^+$, $\Lambda^- =(-\Lambda)\cap
\R^+$. From now on we may assume $\Lambda \subset (0, \infty)$.
Observe that traditionally in the definition of substantial
family, one writes $\sum_k(b_k-a_k)^2a_k^{-2}=\infty$, while for
our purposes it is more convenient to use
$\sum_k(b_k-a_k)^2b_k^{-2}=\infty$. However, these definitions are
equivalent, since as remarked in \cite{Ko}, the divergence of the
first series implies the divergence of the other.

According to the definitions, $D_{BM}(\Lambda)=+\infty$ means that
for every $D>0$ there exists a substantial family of intervals for
$D$. In order to quantify this we introduce the following
definition.

\begin{defnn}
Suppose $\Lambda\subset(0,\infty).$ If $\Psi (s)$ is an increasing
function in $(0, +\infty)$, we call a family of disjoint intervals
$I_k=(a_k, b_k)\subset(0,\infty)$ substantial for $\Psi$ if
$$
\text{$\frac{\Lambda (a_k,b_k)}{b_k-a_k} > \Psi(b_k-a_k)$ \; and \; $\sum_k \bigl (
\frac{b_k - a_k}{b_k}\bigr )^2 = + \infty$}\,.
$$
\end{defnn}

By a diagonal procedure it is immediate to prove:

\begin{lem}\label{l31}
If $\Lambda \subset (0, \infty)$ and $D_{BM} (\Lambda)=+\infty$,
there exists an increasing function $\Psi (s) \nearrow + \infty$
which admits a substantial family of intervals.
\end{lem}

This lemma gives all we will use on $\Lambda$. To prove the
implication (a) $\Rightarrow$ (b), for a given $\Lambda$ with
$D_{BM}(\Lambda)=\infty,$ we have to find a function
$\sigma(y)\nearrow\infty$ such that $\Lambda$ is a uniqueness set
for $B_\sigma.$ This is the subject of the following

\begin{thm}\label{t}
Assume $\Lambda\subset(0,\infty)$, and a family of disjoint
intervals $ I_k=(a_k,b_k)$ is substantial for some function
$\Psi(s)\nearrow\infty$. Then $\Lambda$ is a uniqueness set for
$B_\sigma$, whenever $\sigma$ satisfies:
\begin{enumerate}
\item[(a)] For all $x\in\R$
\begin{equation}\label{12}
\sigma(x)\leq\frac{1}{2e}\Psi(\frac{x}{2e}).\
\end{equation}
\item[(b)] There exists a sequence of integers $n_j\to\infty$ such that
\begin{equation}\label{13}
\frac{1}{\sigma(2b_{n_j})}\sum_{k=1}^{n_j}
\Psi(b_k-a_k)\left(\frac{b_k-a_k}{b_k}\right)^2\to\infty,\
 \ j\to\infty.
\end{equation}

\end{enumerate}
\end{thm}

It is clear that for any function $\Psi(s)\nearrow\infty$ there is
a function $\sigma(s)\nearrow\infty$ which satisfies both
assumptions (\ref{12}) and (\ref{13}). Hence, the implication (a)
$\Rightarrow$ (b) of Theorem~2 is a consequence of Theorem
\ref{t}.

\newpage

\noindent 3.3. The proof of Theorem \ref{t} will be a consequence
of a series of lemmas.

Denote by $B_\sigma^0$ the subclass of $B_\sigma$ of functions
\begin{equation}\label{bss}
|F(x+iy)|\leq e^{|y|\sigma(|y|)},\ x+iy\in\mathbb{C}.
\end{equation}

\begin{lem}\label{l34}
Assume a function $F\in B_\sigma^0$  has $n$ zeros on some
interval $[a,b]$. Then
\begin{equation}\label{14}
|F(x)|\leq (b-a)^n\min_{y>0}\frac{e^{y\sigma(y)}}{y^n}, \ \
\mbox{for every } x\in[a,b].
\end{equation}
\end{lem}

\begin{proof}
Indeed, set $Q(x)=\prod_{k=1}^n(x-x_k)$ where $x_k$ are the zeros
of $F$ on $(a,b)$. The following interpolation formula is well
known: for every $x\in [a,b]$ there is a number $t\in (a,b)$ such
that
$$F(x)=\sum_{k=1}^n\frac{F(x_k)Q(x)}{Q'(x_k)(x-x_k)}+
\frac{F^{(n)}(t)Q(x)}{n!}=\frac{F^{(n)}(t)Q(x)}{n!}.$$ By Cauchy's
inequality, we have $$|F(x)|\leq |Q(x)| \max_{a\leq t\leq b}
\frac{|F^{(n)}(t)|}{n!}\leq (b-a)^n\max_{a\leq t\leq
b}\min_{R>0}\frac{\max_{0\leq\theta<2\pi}|F(t+Re^{i\theta})|}{R^n}$$$$\leq
(b-a)^n\min_{y>0}\frac{e^{y\sigma(y)}}{y^n}.
$$
\end{proof}

\begin{lem}\label{l35}
Suppose $\Psi$ and $\sigma$ satisfy (\ref{12}). If $F\in
B_\sigma^0$ has $n\geq (b-a)\Psi(b-a)$ zeros on an interval
$(a,b), a>0,$ then
\begin{equation}\label{15}
\int_a^b\frac{\log|F(x)|}{x^2}\,dx \leq -(\log
2)(\frac{b-a}{b})^2\Psi(b-a).
\end{equation}
\end{lem}

\begin{proof}
 Let $y_n$ be the number defined by $y_n\sigma(y_n)=n$. Then, by
(\ref{14}),
$$|F(x)|\leq (b-a)^n\frac{e^{y_n\sigma(y_n)}}{y_n^n}= (b-a)^n e^{n-n\log y_n}, \ \
\mbox{for every } x\in[a,b].$$ Hence,
$$\int_a^b\frac{\log|F(x)|}{x^2}\, dx \leq
\frac{b-a}{ab}\left(n\log(b-a)+n-n\log
y_n\right)=-\frac{b-a}{ab}\left(n\log\frac{y_n}{e(b-a)}\right).
$$
Since $n\geq (b-a)\Psi(b-a)$,  it follows from (\ref{12}) and the
definition of $y_n$ that $y_n\geq 2 e(b-a)$, and this gives
(\ref{15}).
\end{proof}

\newpage

Lemma 3.5. relates the size of the logarithmic integral to the
number of zeros. One can obtain a slightly better estimate by
using Jensen's formula for ellipses (see \cite{Ko} for the case of
Bernstein classes). On the other hand, the next lemma establishes
that the logarithmic integral of $F\in B_\sigma$ cannot be too
large negative.

\begin{lem}\label{l36}
If $F\in B_\sigma^0$, $F\not\equiv 0$ then
\begin{equation}\label{16}
  \int_{1\leq |x|\leq R} \frac{\log |F(x)|}{x^2} dx \geq -\frac{16}{3}\, \sigma
  (2R\,)+ 0(1),\, \ R\to\infty.
\end{equation}
\end{lem}

\begin{proof}
We use Carleman's formula \cite[ch.~5]{Le}; for $R>1$:
$$
\int_1^R \Bigl ( \frac{1}{x^2}- \frac{1}{R^2} \Bigr ) \log \,
|F(x)\,F(-x)|\, dx + \frac{2}{R} \int_0^\pi \log
\,|F(Re^{i\theta})|\sin \theta\, d\,\theta \geq C
$$
where $C$ depends only on $F$. By (\ref{bss}), $\log |F(x)|\leq
0$. Using that $x^{-2} \leq \frac 43 (x^{-2}- R^{-2})$ for
$|x|\leq R/2$ we get
$$
\frac 34 \int_{1\leq |x|\leq R/2} \frac{\log |F(x)|}{x^2} \geq C -
4 \sigma (R) \,.
$$
\end{proof}
In fact the proof of the lemma shows that for $F$ satisfying
$|F(z)|=0 \bigl ( e^{\sigma (|z|)|z|} \bigr )$, $\log^- |F(x)|$
controls $\log^+ |F(x)|$.

With Lemmas 3.5. and 3.6. we can now complete the proof of Theorem
3.3.  Suppose $G\not\equiv 0 \in B_\sigma$ vanishes on $\Lambda$.
Set $F=G/C$ where $C$ is a constant such that $F\in B_\sigma^0.$
Then Lemma 3.5. applies to every $I_k=[a_k,b_k]$, so adding
(\ref{15}) in $k\leq n$
$$
\int_1^{b_n} \frac{\log\, |F(x)|}{x^2} dx \,\leq \,- C_1
\sum_{k=1}^n \Psi\, (b_k-a_k) \,\bigl (\frac{b_k-a_k}{b_k}\bigr
)^2,
$$
which combined which (\ref{16}) gives
$$
\sigma (2b_n)- C_2 \sum_{k=1}^n \Psi (b_k-a_k) \bigl (
\frac{b_k-a_k}{b_k}\bigr )^2 \geq C_3,
$$
with constants $C_1$, $C_2$ and $C_3$ independent of $n$. This is
in contradiction with hypothesis (\ref{13}). Hence $F\equiv 0$,
and $\Lambda$ is a uniqueness set for $B_\sigma$.

\newpage

\section{Pairs of generators}
\begin{defn} We say that $\Lambda$ admits a pair of generators if there
 exist two $L^1-$functions $\varphi_1$  and $\varphi_2$ such that
 all linear combinations of  $\varphi_1(t-\lambda_1)$ and
 $\varphi_2(t-\lambda_2),\lambda_1,\lambda_2\in \Lambda,$
 are dense in $L^1 (\R)$.
 \end{defn}

Suppose a sequence $\Lambda$ has a finite spectral radius. Then by
Theorem 1, the $\Lambda-$translations of one function cannot be
dense in $L^1(\R)$. However, $\Lambda$ may admit a pair of
generators. Observe that all perturbations of the integers (2)
have spectral radius $\pi$ (that is, $D_{BM}(\Lambda)=1$). If the
$a_n$ in (2) are exponentially small
 then $\Lambda$ does admit a pair of generators, as shown in the
 following

\begin{thm}\label{pair}
Suppose a real sequence $a_n, n\in \mathbb{Z},$ satisfies:
\begin{equation}\label{small}0\ne |a_n|\leq C r^{|n|},\ \ n\in \mathbb{Z},\end{equation} where $C>0$
and $0<r<1$ are some constants. Then the sequence
$\Lambda=\{n+a_n\}_{n\in \mathbb{Z}}$ admits a pair of generators.
\end{thm}

It is easy to check that the set of integers itself does not admit
a pair (nor any finite number) of generators. Thus the situation
with two generators in $L^1(\R)$ is somewhat similar to the
situation with one generator in $L^2(\R)$, where the set of
integers does not admit generator while  every small perturbation
of it does.

We shall need two auxiliary lemmas.

\begin{lem}\label{sis}
Suppose $a_n$ satisfy (\ref{small}), $\rho<\pi,$ and a function
$g\in B_\rho.$ If there is a constant $C>0$ such that
$|g(n+a_n)|\leq C|a_n|$ for all $n\in\Z$ then $g\equiv
0$.\end{lem}

\begin{proof}
Let $b>0$ satisfy $\rho+b<\pi$. By (\ref{sinc}), the function
$g(z)\,\sinc \, bz$ belongs to $B_{\rho+b}\cap L^2(\R)$ (that is
it belongs to the Paley-Wiener space $PW_{\rho+b}$). Let us
estimate the values of $g$ at the integer points:
\begin{equation}\label{small1}|g(n)|\leq |g(n+a_n)|+|g(n+a_n)-g(n)|\leq (C+\Vert
g^{'}\Vert_\infty)|a_n|= K|a_n|,
\end{equation}
where $K<\infty$ since by Bernstein's inequality $\Vert
g^{'}\Vert_\infty\leq (2\rho+b)\Vert g\Vert_\infty<\infty$.

Set $G=\hat g$. Then $G \in L^2(\R)$ and is concentrated on
$[-\rho-b,\rho+b]$. Write
$$G(x)=\frac{1}{2\pi}\sum_{n=-\infty}^\infty g(-n)e^{inx},$$
where the Fourier series converges to $G$ in $L^2$ norm. One can
easily verify that by (\ref{small}) and (\ref{small1}) the Fourier
series admits analytic continuation into some strip in the complex
plane. However, since $G(x)=0$ on $(\rho+b,\pi)$, the Fourier
series is identically zero, and we conclude that $G(x)=0$ a.e.
\end{proof}

 To formulate the next lemma, we introduce two auxiliary functions
$$\varphi_{1}(t)=\sinc^2at\sum_{k=-\infty}^\infty e^{-|k|+2\pi i k t}, \ \ \varphi_{2}(t)=e^{-i\pi t}\varphi_{1}(t),$$
where $a$ is a constant to be chosen later.

\begin{lem}\label{sist}
Suppose $0<a<\pi/2,$ and $a_n$ satisfy (\ref{small}). Then the set
$\Lambda=\{n+a_n\}$ is a uniqueness set for each class
$L^\infty\star\check\varphi_j, j=1,2.$
\end{lem}

\begin{proof} We use an argument
similar to the one in the proof of the main result of \cite{OU}.
Clearly, $\varphi_{1}\in L^1 (\R),$ for all values $a$. Since
$\sinc\, at$ has as Fourier transform the characteristic function
$\chi_a (x)$ of $[-a,a]$, we get
$$\hat\varphi_{1}=\chi_a\star\chi_a\star\sum_{k=-\infty}^\infty e^{-|k|}\delta_{2\pi k},$$
where $\delta_{c}$ is the unite measure concentrated at the point
$c$. This shows that the support of  the Fourier transforms of
$\varphi_j$ is as follows:
\begin{equation}\label{support}
\mbox{supp
}\hat\varphi_1=\bigcup_{k=-\infty}^\infty[-2a+2k\pi,2a+2k\pi]=[-2a,2a]+2\pi\Z,
\ \mbox{supp }\hat\varphi_2=\mbox{supp }\hat\varphi_1+\pi.
\end{equation}
Observe that both $\hat\varphi_1$ and $\hat\varphi_2$ are strictly
positive on their support.

We have to verify that $f\in L^\infty(\R)$ and
$f\star\check\varphi_{j}(n+a_n)=0, $ for all $n\in\Z$, imply
$f\ast\check\varphi_{j}= 0$. However, it suffices to check this
for $j=1$ only. Indeed, since $\varphi_2(t)=e^{-\pi i
t}\varphi_1(t)$, we have $f\star\check\varphi_2(t)=e^{\pi i
t}[(fe^{-\pi i s})\star\check\varphi_1](t)$, and the implication
for $j=2$ follows from the one for $j=1$. Notice too that
$\varphi_1$ is even, that is, $\check\varphi_1=\varphi_1$.

Set $$g_j(t)=\sum_{k=-\infty}^\infty
e^{-|k|}k^j\int_{-\infty}^\infty \sinc^2 a(t-s) e^{-2\pi i
ks}f(s)\,ds, \ j=0, 1, 2, ...$$ so that
$$f\star\varphi_1
(t)=\sum_{k=-\infty}^\infty e^{-|k|+2\pi ikt}\int_{-\infty}^\infty
\sinc^2 a(t-s) e^{-2\pi i ks}f(s)\,ds=
\sum_{j=0}^\infty\frac{(2\pi i t)^j}{j!}g_j(t).$$ To prove the
lemma, we show by induction that $f\star\varphi_1(n+a_n)=0,
n\in\Z,$ implies
\begin{equation}\label{ind}
g_j\equiv 0, \ j=0, 1, 2, ...
\end{equation}

By (\ref{sinc}) we see that $g_j\in B_{2a}$ for all $j=0, 1, ...$
Since
$$0=f\ast\check\varphi_1(n+a_n)=\sum_{k=-\infty}^\infty
e^{-|k|+2\pi i k
a_n}\int_{-\infty}^\infty\sinc^2\,a(n+a_n-s)e^{-2\pi i
ks}f(s)\,ds,$$ we obtain for each $n\in\Z$
$$|g_0(n+a_n)|=\left|\sum_{k=-\infty}^\infty e^{-|k|}\left(e^{2\pi i k a_n}-1\right)
\int_{-\infty}^\infty\sinc^2\,a(n+a_n-s)e^{-2\pi i
ks}f(s)\,ds\right|\leq C_0 |a_n|,$$ where
$$C_0= \Vert f\Vert_\infty\Vert\sinc^2 at\Vert_1\sum_{k=-\infty}^\infty 2\pi |k| e^{-|k|}<\infty.$$
Hence, by Lemma \ref{sis}, we see that (\ref{ind}) holds for
$j=0$.

Suppose (\ref{ind}) is true for $0\leq n<l.$ Clearly, for each
integer $l$ there is a constant $K_l$ such that the inequality
$$|e^{i\alpha}-\sum_{j=0}^{l}\frac{(i\alpha)^j}{j!}|\leq K_l
|\alpha|^{l+1}$$  holds for every real $\alpha$. Applying it with
$\alpha=2\pi k a_n$, we obtain:
$$\frac{(2\pi a_n)^l}{l!}|g_l(n+a_n)|=\left|f\star\check\varphi_1(n+a_n)-\sum_{k=0}^{l}\frac{(2\pi i a_n)^j}{j!}g_j(n+a_n)\right|
=$$
$$\left|\sum_{k=-\infty}^\infty e^{-|k|}\left(e^{2\pi i ka_n}-\sum_{j=0}^{l}\frac{(2\pi i ka_n)^j}{j!}\right)\int_{-\infty}^\infty
\sinc^2 a(n+a_n-s)e^{-2\pi i ks}f(s)ds\right|$$
$$\leq C_l |a_n|^{l+1},\ \mbox{where }\ C_l=K_l(2\pi)^{l+1}\sum_{k=-\infty}^\infty
e^{-|k|}|k|^{l+1}\Vert f\Vert_\infty\Vert \sinc^2
at\Vert_1<\infty.$$Hence, $|g_l(n+a_n|\leq C_l|a_n|, n \in\Z$
(note that here is where the assumption $a_n\ne 0$ is used). Lemma
\ref{sis} gives $g_l\equiv 0$, so that (\ref{ind}) is true, which
proves the lemma.
\end{proof}

Let us now turn to the proof of Theorem \ref{pair}.
 We shall now  simply check that for every
$\pi/4<a<\pi/2$ the functions $\varphi_1$ and $\varphi_2$ form a
pair of generators for every sequence $\Lambda=\{n+a_n\}$
satisfying the assumptions of  Theorem \ref{pair}. Assume a
function $f\in L^\infty(\R)$ satisfies $f\star\check
\varphi_j(n+a_n)=0,j=1,2,$ for all $n\in\Z$. To prove Theorem
\ref{pair} we have to show that $f=0$ a.e. By  Lemma \ref{sist},
$f\star\check\varphi_j= 0$, $j=1,2.$ Hence,
$f\star(\check\varphi_1+\check\varphi_2)=0$. However, as it easily
follows from (\ref{support}), the function
$\hat\varphi_1+\hat\varphi_2$ is everywhere positive for
$a>\pi/4.$ We conclude, by Wiener's theorem, that $f=0$ a.e.


\begin{thebibliography}{NSV}

\bibitem[AO]{AO} A.\ Atzmon \& A.\ Olevskii, \emph{Completeness of integer
translates in Function Spaces an $\R$}, J. Approximation Theory 83
\textbf{3}(1996)
%
\bibitem[CDH]{CDH} O.\ Christensen, B.\ Deng \& C.\ Heil, \emph{Density of Gabor
frames}, Applied and Computational Harmonic Analysis, \textbf{7},
292-304(1999).
%
\bibitem[F]{F} B.\ Faxen, \emph{On approximation by equidistant
translates}, Pre-print, 1996. Royal Inst. of Technology,
Stockholm, Sweden
%
\bibitem[Ko]{Ko} P.\ Koosis, \emph{The Logarithmic integral}, Cambridge
Univ.\ Press (1992) Vol.~I--II.
%
\bibitem[La]{La} H.\ J.\ Landau, \emph{A sparse sequence of exponentials
closed on large sets}, Bull.\ Amer.\ Math.\ Soc.\ \textbf{70}
(1964), 566--569.
%
\bibitem[Le]{Le} B.\ Ya.\ Levin, \emph{Lectures on entire functions},
Amer.\ Math.\ Soc.\ Providence, (1996).
%
\bibitem[N]{N} N.\ Nikolskii, \emph{Selected problemes of wighted
approximation and spectral analysis}, Trudy Math.\ Inst.\ Steklov
\textbf{120} (1974); English transl: Proc.\ Steklov Math.\ Inst.\
\textbf{120} (1974); Arts Providence, 1976.
%
\bibitem[O]{O} A.\ Olevskii, \emph{Completeness in $L^2 (\R)$ of almost
integer translates}, C.R. Acad.\ Sci.\ Paris, \textbf{324},
S\'{e}rie~1 (1997), 987-991.
%
\bibitem[OU]{OU} A.\ Olevskii, \& A.\ Ulanovskii, \emph{Almost Integer
Translates. Do Nice Generators Exist?} To be published in J.\ of
Fourier Anal.\ and Appl.
%
\bibitem[S]{S} A.M.\ Sedletskii{Approximation by Translations of Functions on the Line} Teor. Pribl. Fuk.
in: Proc. of International conf. Kiev, 31 May - 5 june, 1983. M.
1987 (Russian).
%
\bibitem[U]{U}  A.\ Ulanovskii, \emph{Sparse sequences of fucntions closed on large sets in $\R^n$} J.
London Math. Soc (2), 63 (2001), 428-440.
%
\bibitem[Z]{Z}  R.A.\ Zalik, \emph{On approximation by shifts and a theorem of
Wiener}, Trans. Amer. Math. Soc. \textbf{243}(1978), 299-308.
\end{thebibliography}
\end{document}